\documentclass[fullpage,twoside]{article}
\usepackage{fancyhdr,graphicx}
\usepackage{amsfonts}
\usepackage{booktabs}
\usepackage{array}
\usepackage{amsmath}
\usepackage{graphicx}
\usepackage{amsmath,bm}
\begin{document}
	\title{\vspace{-1cm} An interactive fuzzy goal programming algorithm to solve decentralized bi-level multiobjective fractional programming problem}
	\author{Hasan Dalman \\
		{\small\em  Deparment of Mathematics Engineering, Yildiz Technical University, Esenler 34210, Turkey}\\
		{\small\em hsandalman@gmail.com}}
	\date{}
	\maketitle
	\maketitle \thispagestyle{fancy} \label{first}
	\vspace{-5mm} \pagestyle{plain}
	\newsavebox{\mytable}

\begin{abstract}
	%\vskip 0.15mm\noindent {\bf Abstract:}
This paper proposes a fuzzy goal programming based on Taylor series for solving decentralized bi-level multiobjective fractional programming (DBLMOFP) problem. In the proposed approach, all of the membership functions are associated with the fuzzy goals of each objective at the both levels and also the fractional membership functions are converted to linear functions using the Taylor series approach. Then a fuzzy goal programming is proposed to reach the highest degree of each of the membership goals by taking the most satisfactory solution for all decision makers at the both levels. Finally, a numerical example is presented to illustrate the effectiveness of the proposed approach. 


\end{abstract}
\section{Introduction}

Multiple level programming problems are usually faced with the much hierarchical structure of large organizations such as government offices, profit or non-profit organizations, manufacturing plants, logistic companies, etc. Solution procedures obviously assign each decision maker a unique objective, a set of decision variables and a set of general constraints that affect all decision makers. Each unit independently investigates itself interest but is affected by the actions of other units.

Multilevel programming was first proposed by Bracken and McGill [1] to express a decentralized noncooperative decision system with one decision maker and multiple decision makers of the same situation in 1973. The Stackelberg method has been generally employed to multilevel programming problems. It concerns many applications in real life such as strategic planning (Bracken and McGill, [2]), resource allocation (Aiyoshi and Shimizu,[3]), and water management (Anandalingam and Apprey, [4]). In 1990, Ben-Ayed and Blair [5] proved that multilevel programming is an nonlinear programming (NP)-hard problem. In order to establish the model mathematically, many methods and algorithms have been presented like extreme point algorithm (Candler and Towersley, [6]), k.th best algorithm (Bialas and Karwan, [7]), branch and bound algorithm (Bard and Falk, [8]), descent method (Savard and Gauvin, [9]). and genetic algorithm (Liu, [10]). However, Gao and Liu [11] suggested a novel fuzzy multilevel programming model and defined a Stackelberg-Nash equilibrium. But, these traditional methods are based on Karush-Kuhn-Tucker conditions and/or penalty functions [12] and besides the Stackelberg method does not produce Pareto optimality because of its non-cooperative nature [13].

Therefore, in a hierarchical decision-making, it has been realized that each decision maker should have a motivation to participate with the other, and a minimum level of fulfillment of the decision maker at a lower level must be remarked for the overall profit of the organization. The use of membership function of fuzzy set theory to multilevel programming problems for satisfactory solutions was first introduced by Lai [12] in 1996. Then Lai's solution procedure was extended by Shih et al. [14, 15]. In addition, fuzzy programming methods were applied by many authors for solving multiple level linear programming problems [16, 17], bi-level quadratic fractional programming problem [13, 17, 18], bi-level nonconvex programming problems with fuzzy parameters [18], decentralized bi-level linear programming problems [19, 20], so far.

Baky [21] employed two fuzzy goal programming method to reach the optimal solutions of multi-level multi-objective linear programming problems. Arora and Gupta [22] joined an interactive fuzzy goal programming algorithm with the theory of dynamic programming for solving the bi-level programming problems. Wang et al. [23] developed an idea to deal with a bi-level multi followers programming problem. The distance function, which was introduced by Yu [24], has been extensively exercised to solve multiobjective programming problem to determine a compromise solution. Based on the theory of distance function, Moitra and Pal [25] applied a fuzzy goal programming procedure and reached a satisfactory balance through minimizing the deviations of the leader and follower as much as possible for bi-level programming. Baky and Abo-Sinna [26] proposed a fuzzy TOPSIS algorithm, which simultaneously minimized a distance function from an idea point and maximized another distance function from a nadir point, to solve bi-level multiobjective decision-making problems. However, some interesting interactive fuzzy decision-making algorithms have widely been employed to bi-level and multi-level programming problems [27, 28, 29, 30, 32, 33] can get preferred result.

Abo-Sinna and Baky [34] suggested a fuzzy goal programming approach to decentralized bi-level multi-objective linear fractional programming (DBLMOLP) problems. Baky [34] applied the fuzzy goal programming method of Pal et al. [35] for linearization of membership functions. Toksari and Bilim [36] presented an interactive fuzzy goal programming method based on Jacobian matrix for solving bi-level multiobjective fractional programming problem. They extended the fuzzy goal programming method proposed by Mohamed [37] to reach the optimal solutions of DBLMOFP problems.

In this paper, a fuzzy goal programming algorithm based on Taylor series approach is proposed to solve decentralized bi-level  multiobjective fractional programming (DBLMOFP) problem. Objective and constraint functions of decision makers are fractional and linear functions, respectively. Until now, many solution procedure are introduced to linearize the fractional and/or nonlinear functions in literature. In the case of a fractional objective programming, the most popular solution procedures are based on linearization methods (see [34, 35, 36, 39] for details). For this reason, in the construction of the problem, optimal values of the objectives are found by determining individual optimal solutions. Then the fractional membership functions which are joined with each fractional objective of the problem are determined and then membership functions joined with each objective are converted to linear membership functions by employing Taylor series approach. Thus, an interactive fuzzy goal programming algorithm is presented for solving the DBLMOFP problem. Finally, two numerical examples are given to demonstrate the practicability of the proposed algorithm. The remaining of this paper is classified as follows. Section 2 gives a DBLMOLFP problem with its formulation and Section 2 suggests a fuzzy decision making algorithm. A numerical example is given in Section 4 and Section 5 concludes the paper.

\section{Problem Formulation}

 In a decentralized  bi-level multiobjective fractional programming (DBLMOFP) problem, two decision makers are located at two distinct hierarchical levels with multiple objectives independently checking a set of decision variables. In a DBLMOFP, the first level decision maker called the leader affects its decision in full view of the second level decision makers called the follower. Each decision maker attempts to optimize its objective function and is affected by the actions of the other decision makers.

 Let the vector of decision variables $X_{1} \in R^{n_{1} } $and $X_{k>1} \in R^{n_{k>1} } $ be under the control of the first decision maker (DM) and the second decision makers (DMs), respectively  $n_{1} ,n_{2} ,...n_{m} \ge 1.$ So the DBLMOFP problem of minimization type should be formulated as follows ([13,19, 20, 34, 36]);

\noindent The upper level: 
\[\mathop{\min }\limits_{x_{1} } F_{1} \left(x\right)=\mathop{\min }\limits_{x_{1} } \left(f_{11} \left(x\right),f_{12} \left(x\right),...,f_{1p_{1} } \left(x\right)\right),\] 

where $\left(x_{1} ,x_{2} ,...,x_{m} \right)$ solves  
the lower level: 
\[\mathop{\min }\limits_{x_{2} } F_{2} \left(x\right)=\mathop{\min }\limits_{x_{2} } \left(f_{21} \left(x\right),f_{22} \left(x\right),...,f_{2p_{2} } \left(x\right)\right),\] 
\begin{equation} \label{GrindEQ__1_} 
\mathop{\min }\limits_{x_{3} } F_{3} \left(x\right)=\mathop{\min }\limits_{x_{3} } \left(f_{31} \left(x\right),f_{32} \left(x\right),...,f_{3p_{3} } \left(x\right)\right), 
\end{equation} 
\[\begin{array}{l} {.} \\ {.} \\ {.} \end{array}\] 
\[\mathop{\min }\limits_{x_{m} } F_{m} \left(x\right)=\mathop{\min }\limits_{x_{m} } \left(f_{m1} \left(x\right),f_{m2} \left(x\right),...,f_{mp_{m} } \left(x\right)\right),\] 
Subject to
\[x\in S=\left\{\left. x\in {\rm R}^{n} \right|A_{1} x_{1} +A_{2} x_{2} +....A_{2} x_{m} \left(\begin{array}{l} {\le } \\ {=} \\ {\ge } \end{array}\right)b,x_{m} \ge 0,b\in R^{s} \right\}\ne \emptyset \] 

where $f_{ij} \left(x\right)=\frac{c_{1}^{ij} x_{1} +c_{2}^{ij} x_{2} +...+c_{m}^{ij} x_{m} +\alpha _{}^{ij} }{d_{1}^{ij} x_{1} +d_{2}^{ij} x_{2} +...+d_{m}^{ij} x_{m} +\beta _{}^{ij} } ,i=1,2,...,m;j=1,2,...,p_{m} $
where $p$ is the number of the objective functions for decision makers at the lower level and $s$ is the number of constraints.$A_{i} $ is an $s\times p$ constant matrix. $S$ is a non-empty, convex and compact set in ${\rm R}^{n} $ ; $d_{1}^{ij} x_{1} +d_{2}^{ij} x_{2} +...+d_{m}^{ij} x_{m} +\beta _{}^{ij} $is greater than zero.

\subsection{ Construction of fractional membership function}

In multiple objective programming, if an imprecise aspiration level is injected to each of the objectives then these fuzzy objectives are expressed as fuzzy goals. According to Pal and Moitra [25], Li and Hu [42] and Bellman and Zadeh [43], the decision is often formulated as follows:

\begin{equation}
\begin{array}{l}
Find\,\,x\,\\
so\,as\,to\,satisfy\,\\
\,\,\,\,\,\,\,\,\,\,{f_{ij}}\left( x \right)\left( \begin{array}{l}
\le \\
\cong \\
\ge 
\end{array} \right)f_{ij}^*,\,\,i = 1,2,...m;\,j = 1,2,...,p_m^{}\\
subject\,\,to\,\\
\,\,\,\,\,\,\,\,\,\,\,x \in S
\end{array}
\end{equation}
where $f_{ij}^{*} $is the perspective goal value for the objective function $f_{ij} \left(x\right).$ $\left(\lower3pt\hbox{\rlap{$\scriptscriptstyle\sim$}}\prec ,\cong ,\lower3pt\hbox{\rlap{$\scriptscriptstyle\sim$}}\succ \right)$ represent different fuzzy relations.  Let $\left(f_{ij}^{*} ,f_{ij}^{\max } \right)$ be the tolerant interval selected to the $ij.th$ objective $f_{ij} \left(x\right).$Then the fuzzy goals are $f_{ij} \left(x\right)\lower3pt\hbox{\rlap{$\scriptscriptstyle\sim$}}\prec f_{ij}^{*} $for the minimization objective where $\lower3pt\hbox{\rlap{$\scriptscriptstyle\sim$}}\prec $ represent the fuzzified inequalities. Thus, the memberhip function is determined as
\begin{equation} \label{GrindEQ__3_} 
\mu _{ij} \left(f_{ij} \left(x\right)\right)=\left\{\begin{array}{l} {1} \\ {\frac{f_{ij}^{\max } -f_{ij} \left(x\right)}{f_{ij}^{\max } -f_{ij}^{*} } } \\ {0} \end{array}\right. ,  \begin{array}{l} {if} \\ {if} \\ {if} \end{array} \begin{array}{l} {f_{ij} \left(x\right)\le f_{ij}^{*} } \\ {f_{ij}^{*} \le f_{ij} \left(x\right)\le f_{j}^{\max } } \\ {f_{ij} \left(x\right)\ge f_{ij}^{\max } } \end{array} 
\end{equation} 
where $f_{ij}^{*} $ is called an ideal value and $f_{ij}^{\max } $ is tolerance limit for$f_{ij}^{} .$

\noindent Similarly, let $\left(f_{ij}^{\min } ,f_{ij}^{*} \right)$ be the tolerant interval selected to the $ij.th$ objective $f_{ij} \left(x\right).$ the fuzzy goals are $f_{ij} \left(x\right)\lower3pt\hbox{\rlap{$\scriptscriptstyle\sim$}}\succ f_{ij}^{*} $for the minimization objective where $\lower3pt\hbox{\rlap{$\scriptscriptstyle\sim$}}\succ $ represent the fuzzified inequalities. Thus, the memberhip function is determined as
\begin{equation} \label{GrindEQ__4_} 
\mu _{ij} \left(f_{ij} \left(x\right)\right)=\left\{\begin{array}{l} {1} \\ {\frac{f_{ij} \left(x\right)-f_{ij}^{\min } }{f_{ij}^{*} -f_{ij}^{\min } } } \\ {0} \end{array}\right. ,  \begin{array}{l} {if} \\ {if} \\ {if} \end{array} \begin{array}{l} {f_{ij} \left(x\right)\ge f_{ij}^{*} } \\ {f_{ij}^{\min } \le f_{ij} \left(x\right)\le f_{j}^{*} } \\ {f_{ij} \left(x\right)\ge f_{ij}^{\min } } \end{array} 
\end{equation} 
In multi-objective decision making, the tolerant interval of the objective may be unknown. Hence, the payoff table is employed to determine it. 

\noindent Here, to facilitate computation, we used the same method with the paper [28] and [42] to determine membership functions. Let us $\mu _{ij} \left(f_{ij} \left(x\right)\right),i=1,2,...m;j=1,2,...,p_{m}^{} $to define the fuzzy goals of the leader and the follower, respectively. Then the fuzzy goals are $f_{ij} \left(x\right)\lower3pt\hbox{\rlap{$\scriptscriptstyle\sim$}}\prec f_{ij}^{*} $for the minimization objective in model {(1)}

\noindent For the sake of simplicity, we suppose that $f_{ij}^{*} $and $f_{ij}^{\max } $are the ideal value and the tolerance limit of the above fractional problem, respectively. For instance, 
\begin{equation} \label{GrindEQ__5_} 
f_{ij}^{\max } =\max _{x\in S} f_{ij} \left(x\right),i=1,2,...m;j=1,2,...,p_{m}^{}  
\end{equation} 
and 
\begin{equation} \label{GrindEQ__6_} 
f_{ij}^{*} =\min _{x\in S} f_{ij} \left(x\right),i=1,2,...m;j=1,2,...,p_{m}^{}  
\end{equation} 
Then the payoff table is determined as in [42].

\subsection{ Linearization fractional membership function using the Taylor series approach}

Here, membership functions associated with each objective are linearized by using Taylor series approach, at first. Then, linear membership functions of fractional membership functions associated with each objective are determined. The suggested procedures for fractional objectives can be continued as follows:

\noindent Obtain $\tilde{x}_{i}^{*} =\left(\tilde{x}_{i1}^{*} ,\tilde{x}_{i2}^{*} ,...,\tilde{x}_{ip_{i} }^{*} \right)$  which is the value that is used to maximize the ij-th membership function $\mu _{ij} \left(f_{ij} \left(x\right)\right)$ associated with ij-th objective $f_{ij} \left(x\right).$

Convert membership functions by using Taylor series approach as follows:
\begin{equation}
{\tilde \mu _{ij}}\left( {{f_{ij}}\left( x \right)} \right) \cong \left[ {{{\left. {\frac{{{\mu _{ij}}\left( {{f_{ij}}\left( {\tilde x_i^*} \right)} \right)}}{{\partial {x_1}}}} \right|}_{\tilde x_i^*}}\left( {{x_1} - \tilde x_{i1}^*} \right) + {{\left. {\frac{{{\mu _{ij}}\left( {{f_{ij}}\left( {\tilde x_i^*} \right)} \right)}}{{\partial {x_2}}}} \right|}_{\tilde x_i^*}}\left( {{x_2} - \tilde x_{i2}^*} \right) + ... + {{\left. {\frac{{{\mu _{ij}}\left( {{f_{ij}}\left( {\tilde x_i^*} \right)} \right)}}{{\partial {x_m}}}} \right|}_{\tilde x_i^*}}\left( {{x_2} - \tilde x_{i{p_m}}^*} \right)} \right]
\end{equation}

Taylor series approach generally gives a relatively good approximation to a differentiable function but only around a given point, and not over the entire domain. Functions $\tilde{\mu }_{ij} \left(f_{ij} \left(x\right)\right)$ approximates the function $\mu _{ij} \left(f_{ij} \left(x\right)\right)$around its optimal solution.

After determined the linear membership function $\tilde{\mu }_{ij} \left(f_{ij} \left(x\right)\right)$ using the above procedures, therefore we will transform the decentralized bi-level multi-objective linear programming (DBLMOP) problem to DBLMOFP.

\section{A Fuzzy Goal Programming Algorithm to DBLMOFP}

The fuzzy goal programming approach to multi-objective programming problems introduced by Mohamed [37]. The fuzzy goal programming approach to multi-objective programming problems is extended by Baky [34] and Toksari and Bilim [36]  to solve the DBLMOP. But the fuzzy goal programming of Mohammed [37] is not every time proper. In the practical approach and design, decision maker generally has an importance requirement for different objectives in  multiple objective programming problems as his or her preference. This can be illustrated by weights [40]. Tiwari et al. [41] employed the weighted sum method to obtain it. 

\begin{equation}
\left\{ {\begin{array}{*{20}{l}}
	{\max  \left( {\sum\limits_{j = 1}^{{p_1}} {w_{1j}^{}} {{\tilde \mu }_{1j}} + \sum\limits_{j = 1}^{{p_2}} {w_{2j}^{}} {{\tilde \mu }_{2j}} + ... + \sum\limits_{j = 1}^{{p_m}} {w_{mj}^{}} {{\tilde \mu }_{mj}}} \right)}\\
	{s.t.\left\{ {\begin{array}{*{20}{l}}
			{{{\tilde \mu }_{ij}}\left( {{f_{ij}}\left( x \right)} \right) \le 1,\,\,\,\,\,\,\,\,\,\,\,\,i = 1,2,...m;j = 1,2,...,p_m^{}}\\
			{{A_1}{x_1} + {A_2}{x_2} + ....{A_2}{x_m}\left( {\begin{array}{*{20}{l}}
					\le \\
					= \\
					\ge 
					\end{array}} \right)b,}\\
			{x \in S}
			\end{array}} \right.}
	\end{array}} \right.
\end{equation}
where (weight) $w_{ij} $ is the relative importance of achieving the aspired levels of the fuzzy goals. Here, $w_{ij}^{} $can be written as $\left(\mathop{\sum }\limits_{j=1}^{p_{1} } w_{1j}^{} +\mathop{\sum }\limits_{j=1}^{p_{2} } w_{2j}^{} +...+\mathop{\sum }\limits_{j=1}^{p_{m} } w_{mj}^{} \right)=1.$ This model can be optimized all objectives as far as probable. But it may give the same results when the weights are mutated. 

In a fuzzy programming, the highest possible value of membership function is always 1. According to the approach of Mohamed [37], the fractional membership functions in {(3)} and {(4)} can be can be determined as the fractional membership goals as follows;
\begin{equation}
{{\mu _{ij}}\left( {{f_{ij}}\left( x \right)} \right) = \frac{{f_{ij}^{\max } - {f_{ij}}\left( x \right)}}{{f_{ij}^{\max } - f_{ij}^*}} + d_{ij}^ -  - d_{ij}^ +  = 1,{\mkern 1mu} {\mkern 1mu} {\mkern 1mu} {\mkern 1mu} {\mkern 1mu} {\mkern 1mu} {\mkern 1mu} {\mkern 1mu} {\mkern 1mu} {\mkern 1mu} {\mkern 1mu} {\mkern 1mu} {\mkern 1mu} {\mkern 1mu} {\mkern 1mu} {\mkern 1mu} {\mkern 1mu} {\mkern 1mu} {\mkern 1mu} {\mkern 1mu} i = 1,2,...m;j = 1,2,...,p_m^{}}
\end{equation}
\begin{equation}
{{\mu _{ij}}\left( {{f_{ij}}\left( x \right)} \right) = \frac{{{f_{ij}}\left( x \right) - f_{ij}^{\min }}}{{f_{ij}^* - f_{ij}^{\min }}} + d_{ij}^ -  - d_{ij}^ +  = 1,{\mkern 1mu} {\mkern 1mu} {\mkern 1mu} {\mkern 1mu} {\mkern 1mu} {\mkern 1mu} {\mkern 1mu} {\mkern 1mu} {\mkern 1mu} {\mkern 1mu} {\mkern 1mu} {\mkern 1mu} {\mkern 1mu} {\mkern 1mu} {\mkern 1mu} {\mkern 1mu} {\mkern 1mu} {\mkern 1mu} {\mkern 1mu} {\mkern 1mu} {\mkern 1mu} i = 1,2,...m;j = 1,2,...,p_m^{}}
\end{equation}
where $d_{ij}^{-} \left(\ge 0\right)$ symbolize the negative deviational variables and $d_{ij}^{+} \left(\ge 0\right)$ symbolize the positive deviational variables. Therefore, the linear membership function based on Taylor series {(7)} can be formulated as the folowing linear functions:
\begin{equation}
{\tilde \mu _{ij}}\left( {{f_{ij}}\left( x \right)} \right) + d_{ij}^ -  - d_{ij}^ +  = 1\,\,\,\,,i = 1,2,...m;j = 1,2,...,p_m
\end{equation} 
Since the maximum possible value of the linear membership goal {(11)} is unity, positive deviation is not possible and/or unnecessary (see, [41]). Here,  the linear function {(11)} can be written as 
\begin{equation} 
{\tilde \mu _{ij}}\left( {{f_{ij}}\left( x \right)} \right) + d_{ij}^ -  = 1,\,\,\,\,i = 1,2,...m;j = 1,2,...,p_m
\end{equation} 
where $d_{ij}^{-} \left(\ge 0\right)$ symbolize the negative deviational variables.

In this paper, we extended the fuzzy goal programming approach, which is suggested by Gupta and Bhattacharjee [44], to obtain the optimal solution for the bi-level problem. They attempted to determine a new weighted fuzzy goal approach for fuzzy goal programming problem by considering only negative deviational variables $d_{ij}^{-} \left(\ge 0\right)$in the goal constraint for the fuzzy multiobjective goal programming problem with aspiration level one, $i=1,2,...m;j=1,2,...,p_{m}^{} .$

\noindent Then this approach is employed to achieve the highest degree of membership for each of the goals by using max-min operator. The weights are also attached to the fuzzy control operator $\lambda $ in the constraints.

\noindent According the idea of Gupta and Bhattacharjee [44] membership goals based on {(7)}, the upper level decision maker problem to obtain the optimal decision vector $x_{i}^{F} =\left(x_{i1}^{F} ,x_{i2}^{F} ,...,x_{ip_{i} }^{F} \right)$ for the problem {(1)} can be formulated as follows:
\begin{equation} 
\left\{ {\begin{array}{*{20}{l}}
	{\max \lambda }\\
	{s.t.\left\{ {\begin{array}{*{20}{l}}
			{w_{1j}^{}\lambda  \le {{\tilde \mu }_{1j}}\left( {{f_{1j}}\left( x \right)} \right),\,\,j = 1,2,...,p_m^{}}\\
			{\lambda  + d_{1j}^ -  = 1,\,\,j = 1,2,...,p_m^{}}\\
			{{A_1}{x_1} + {A_2}{x_2} + ....{A_m}{x_m}\left( {\begin{array}{*{20}{l}}
					\le \\
					= \\
					\ge 
					\end{array}} \right)b,}\\
			\begin{array}{l}
			\lambda  \ge 0,\\
			x \in S
			\end{array}
			\end{array}} \right.}
	\end{array}} \right.
\end{equation} 
where $w_{1j}^{} $is taken as $w_{1j}^{} =\frac{1}{\left|f_{1j}^{\max } -f_{1j}^{*} \right|} ,$  $j=1,2,...,p_{m}^{} $and $\left|f_{1j}^{\max } -f_{1j}^{*} \right|$ is the tolerance which is relatively taken.  

\noindent In the former studies, Baky did not mentioned about the bounds on the maximum negative and positive tolerance values and appropriate method to determine these values in [34]. Also, satisfactory solution of DBLMOFP problem using algorithms found by Baky [34] depends upon the choice of these tolerance values which many times leads to the possibility of rejecting the solution again and again by upper level decision maker and so the solution process spend too much time [33]. Therefore, Lachhwani [33] presented a linear membership function to determine the linear membership functions of decision variables. 

\noindent From the interactive mechanism of {(1)}, the upper level decision variable should be assigned to the lower level decision makers. After determining the optimal solution $x_{i}^{F} =\left(x_{i1}^{F} ,x_{i2}^{F} ,...,x_{ip_{i} }^{F} \right)$ for the upper level of problem {(1)}, we formulated the suggested fuzzy goal programming to solve DBLMOLP problem. It can be constructed as 
\begin{equation} \label{GrindEQ__14_} 
\begin{array}{l} {\max \lambda } \\ {\left\{\begin{array}{l} {w_{ij}^{} \lambda \le \tilde{\mu }_{ij} \left(f_{ij} \left(x\right)\right),i=1,2,...m;j=1,2,...,p_{m}^{} } \\ {\lambda +d_{ij}^{-} =1,i=1,2,...m;j=1,2,...,p_{m}^{} } \\ {A_{1} x_{1} +A_{2} x_{2} +....A_{m} x_{m} \left(\begin{array}{l} {\le } \\ {=} \\ {\ge } \end{array}\right)b,} \\ {x_{1}^{} =x_{i1}^{F} } \\ {\lambda \ge 0,x\in S} \end{array}\right. } \end{array} 
\end{equation} 
Here, $w_{ij}^{} $is considered as $w_{ij}^{} =\frac{1}{\left|f_{ij}^{\max } -f_{ij}^{*} \right|} ,$  $i=1,2,...m;j=1,2,...,p_{m}^{} .$ 

\noindent So the lower level solution can be obtained as $x_{i}^{S} =\left(x_{i1}^{F} ,x_{i2}^{S} ,...,x_{ip_{i} }^{S} \right).$

\noindent 
\subsection{ The suggested fuzzy goal programming algorithm to solve DBLMOLFP}

\noindent \textbf{Step 1} Solve the problem {(1)} as in equation {(5)} and {(16)} by taking single objective function at a time and neglecting all others and obtain the optimal solutions$x_{i}^{*} =\left(x_{i1}^{*} ,x_{i2}^{*} ,...,x_{ip_{i} }^{*} \right)$. 

\noindent \textbf{Step 2} Then determine the ideal value(s)$f_{ij}^{*} ,$$i=1,2,...m;j=1,2,...,p_{i}^{} ,$ tolerans limit(s) $f_{ij}^{\max } ,$$i=1,2,...m;j=1,2,...,p_{i}^{} $ and the relative importance $w_{ij}^{} $ for each objective at each levels. 

\noindent \textbf{Step 3} Construct the fractional membership function {(3)} and/ or {(4)}for each objective at each levels.

\noindent \textbf{Step 4 }Determine $\tilde{x}_{i}^{*} =\left(\tilde{x}_{i1}^{*} ,\tilde{x}_{i2}^{*} ,...,\tilde{x}_{ip_{i} }^{*} \right)$which is the value that is used to maximize the $ij-$ th membership function $\mu _{ij} \left(f_{ij} \left(x\right)\right)$ associated with objective $f_{ij} \left(x\right)$at each levels.

\noindent  \textbf{Step 5 }Then linearize all fractional membership functions using first-order Taylor polynomial series around the decision vector$\tilde{x}_{i}^{*} =\left(\tilde{x}_{i1}^{*} ,\tilde{x}_{i2}^{*} ,...,\tilde{x}_{ip_{i} }^{*} \right)$.

\noindent \textbf{Step 6\textit{ }}Formulate the fuzzy goal programming model {(13)}to obtain the optimal solution for the upper level of {(1)}. Then, solve it to obtain the alternative optimal solution $x_{i}^{F} =\left(x_{i1}^{F} ,x_{i2}^{F} ,...,x_{ip_{i} }^{F} \right)$.

\noindent \textbf{Step 7\textit{ }}Assigning the upper level decision variables to the lower level decision making problem , formulate the fuzzy goal programming model {(14)} and then, solve it to obtain the candidate optimal solution $x_{i}^{S} =\left(x_{i1}^{F} ,x_{i2}^{S} ,...,x_{ip_{i} }^{S} \right).$

\noindent \textbf{Step 8 }If the decision maker for the upper level decision maker is satisfied by the current candidate solution $x_{i}^{S} =\left(x_{i1}^{F} ,x_{i2}^{S} ,...,x_{ip_{i} }^{S} \right)$in Step 9, go to Step 9, else go to Step 10.

\noindent \textbf{Step 9 }The candidate solution $x_{i}^{S} =\left(x_{i1}^{F} ,x_{i2}^{S} ,...,x_{ip_{i} }^{S} \right)$is the optimal solution for the decentralized bi-level multi-objective fractional programming problems (DBLMOFPP).

\noindent \textbf{Step 10 }Modify the tolerance limit $f_{ij}^{\max } ,$ $i=1,2,...m;j=1,2,...,p_{m}^{} $  and the relative importance $w_{ij}^{} ,$$i=1,2,...m;j=1,2,...,p_{m}^{} $ for all objective at all levels and go to Step 3.

\section{Numerical Example}

The suggested interactive fuzzy goal programming approach will be used to a known numerical example. The following bilevel decentralized multiobjective programming problem solved by Baky [34] using interactive fuzzy goal programming method. Now, we will solve this problem using the suggested approach in this paper. 

\textbf{Example 1}

Consider the following problem

Upper level:\textbf{ }
\[\mathop{\min }\limits_{x_{0} } \left(f_{11} =\frac{-x_{0} -4x_{1} +x_{2} +1}{2x_{0} +3x_{1} +x_{2} +2} ,f_{12} =\frac{-2x_{0} +x_{1} +3x_{2} +4}{2x_{0} -x_{1} +x_{2} +5} \right)\] 

Lower level:\textbf{ }

[DM1]    $\mathop{\min }\limits_{x_{1} } \left(f_{21} =\frac{3x_{0} -2x_{1} +2x_{2} }{x_{0} +x_{1} +x_{2} +3} ,f_{22} =\frac{-7x_{0} -2x_{1} +x_{2} +1}{5x_{0} +2x_{1} +x_{2} +1} \right)$

[DM2]    $\mathop{\min }\limits_{x_{2} } \left(f_{31} =\frac{x_{0} +x_{1} +x_{2} -4}{x_{0} -2x_{1} +10x_{2} +6} ,f_{32} =\frac{2x_{0} -x_{1} +x_{2} +4}{-x_{0} +x_{1} +x_{2} +10} \right)$\textbf{}

Subject to
\[\begin{array}{l} {g_{1} =x_{0} +x_{1} +x_{2} \le 5,} \\ {g_{2} =x_{0} +x_{1} -x_{2} \le 2,} \\ {g_{3} =x_{0} +x_{1} +x_{2} \ge 1,} \end{array}    \begin{array}{l} {g_{4} =-x_{0} +x_{1} +x_{2} \le 1,} \\ {g_{5} =x_{0} -x_{1} +x_{2} \le 4,} \\ {g_{6} =x_{0} +2x_{2} \le 4,} \\ {x_{0} ,x_{1} ,x_{2} \ge 0.} \end{array}\] 

(Step 1 and Step 2): Table 2 presents the ideal values, tolerance limits and weights of all the objective functions in both the levels. 

\begin{center} Table 2 The individual minimum and maximum values, the ideal value and tolerance limits and weights

\begin{tabular}{|p{0.7in}|p{0.7in}|p{0.7in}|p{0.7in}|p{0.7in}|p{0.7in}|p{0.7in}|} \hline 
 & $f_{11} $ & $f_{12} $ & $f_{21} $ & $f_{22} $ & $f_{31} $ & $f_{32} $ \\ \hline 
$\max f_{ij} $ & $0.67$ & $1.25$ & $1.353$ & $1$ & ${\rm -0.026}$ & $1.125$ \\ \hline 
$\min f_{ij} $ & $-0.733$  & $0$  & ${\rm -0.50}$  & ${\rm -1.18}$  & $-0.75$  & $0.27$  \\ \hline 
$f_{ij}^{\max } $ & $0.6$ & $1.2$ & $1.3$ & $1$ & ${\rm -0.05}$ & $1.125$ \\ \hline 
$f_{ij}^{*} $ & $-0.7$  & $0$  & ${\rm -0.50}$  & ${\rm -1}$  & $-0.75$  & $0.25$  \\ \hline 
$w_{ij}^{} $ & $0.769$ & ${\rm 0.83}$ & ${\rm 0.56}$ & ${\rm 0.5}$ & $1.43$ & $1.143$ \\ \hline 
\end{tabular}
\end{center}
(Step 3): Thus, the fractional membership functions $\mu _{ij} \left(f_{ij} \left(x\right)\right)$ at each levels are constructed as

\[\begin{array}{l}
{\mu _{11}}\left( {{f_{11}}\left( x \right)} \right) = 0.46 - \frac{{0.77\left( { - {x_0} - 4{x_1} + {x_2} + 1} \right)}}{{2{x_0} + 3{x_1} + {x_2} + 2}},\\
{\mu _{12}}\left( {{f_{12}}\left( x \right)} \right) = 1.00 - \frac{{0.83\left( { - 2{x_0} + {x_1} + 3{x_2} + 4} \right)}}{{2{x_0} - {x_1} + {x_2} + 5}},\\
{\mu _{21}}\left( {{f_{21}}\left( x \right)} \right)0.72 - \frac{{0.56\left( {3{x_0} - 2{x_1} + 2{x_2}} \right)}}{{{x_0} + {x_1} + {x_2} + 3}},\\
{\mu _{22}}\left( {{f_{22}}\left( x \right)} \right) = 0.50 - \frac{{0.50\left( { - 7{x_0} - 2{x_1} + {x_2} + 1} \right)}}{{5{x_0} + 2{x_1} + {x_2} + 1}},\\
{\mu _{31}}\left( {{f_{31}}\left( x \right)} \right) =  - 0.07 - \frac{{1.43\left( {{x_0} + {x_1} + {x_2} - 4} \right)}}{{{x_0} - 2{x_1} + 10{x_2} + 6}},\\
{\mu _{32}}\left( {{f_{32}}\left( x \right)} \right) = 1.125 - \frac{{1\left( {2{x_0} - {x_1} + {x_2} + 4} \right)}}{{ - {x_0} + {x_1} + {x_2} + 10}}.
\end{array}\]
(Step 4): Furthermore, the maximal solution points of fractional membership functions are given in Table 4.

\begin{center}Table 3 The optimal  solutions for all membership functions in both levels under the constraints

\begin{tabular}{|p{0.7in}|p{0.7in}|p{0.7in}|p{0.7in}|p{0.7in}|p{0.7in}|p{0.7in}|} \hline 
 & $\mu _{11} \left(f_{11} \left(x\right)\right)$ & $\mu _{12} \left(f_{12} \left(x\right)\right)$ & $\mu _{21} \left(f_{21} \left(x\right)\right)$ & $\mu _{22} \left(f_{11} \left(x\right)\right)$ & $\mu _{31} \left(f_{31} \left(x\right)\right)$ & $\mu _{32} \left(f_{32} \left(x\right)\right)$ \\ \hline 
$\left(x_{0}^{*} ,x_{1}^{*} ,x_{2}^{*} \right)$ & $\left(0.5,1.5,0\right)$ & $\left(2,0,0\right)$ & $\left(0,1,0\right)$ & $\left(2,0,0\right)$ & $\left({\rm 0,1,0}\right)$ & $\left({\rm 0,1,0}\right)$ \\ \hline 
\end{tabular}
\end{center}

 (Step 5): Then, we employed (from Table 3) to linearize all fractional membership functions using first order Taylor series approach {(7)} around the individual optimal points. Then, we obtained the following linear functions as
\[\left\{ {\begin{array}{*{20}{l}}
	{\begin{array}{*{20}{l}}
		{{{\tilde \mu }_{11}}\left( {{f_{11}}\left( x \right)} \right) = 0.773 - 0.049{x_0} + 0.185{x_1} - 0.178{x_{{2_2}}},}\\
		{{{\tilde \mu }_{12}}\left( {{f_{12}}\left( x \right)} \right) = 0.630 + 0.185{x_0} - 0.093{x_1} - 0.278{x_2},}
		\end{array}}\\
	{{{\tilde \mu }_{21}}\left( {{f_{21}}\left( x \right)} \right) = 0.792 - 0.486{x_0} + 0.208{x_1} - 0.347{x_2},}\\
	{{{\tilde \mu }_{22}}\left( {{f_{22}}\left( x \right)} \right) = 0.992 + 0.050{x_0} - 0.017{x_1} - 0.099{x_2},}\\
	{{{\tilde \mu }_{31}}\left( {{f_{31}}\left( x \right)} \right) = 0.821 - 0.625{x_0} + 0.179{x_1} - 3.036{x_2},}\\
	{{{\tilde \mu }_{32}}\left( {{f_{32}}\left( x \right)} \right) = 0.737 - 0.207{x_0} + 0.116{x_1} - 0.066{x_2},}
	\end{array}} \right.\]
(Step 6): Consequently, the considered fuzzy goal programming model for the upper level decision maker is constructed as follows based on model {(13)}:
The alternative optimal solutions are determined as
\[x_{}^{F} =\left(x_{0}^{F} =1.25,x_{1}^{F} =0.31,x_{2}^{F} =0\right)\] 

(Step 7): The suggested fuzzy goal programming model to DBLMOFPP continues as follows based on {(14)};

\[\left\{ \begin{array}{l}
\max \lambda \\
s.t.\left\{ \begin{array}{l}
0.769\lambda  \le 0.773 - 0.049{x_0} + 0.185{x_1} - 0.178{x_2},\\
0.83\lambda  \le 0.630 + 0.185{x_0} - 0.093{x_1} - 0.278{x_2},\\
0.56\lambda  \le 0.792 - 0.486{x_0} + 0.208{x_1} - 0.347{x_2},\\
0.5\lambda  \le 0.992 + 0.050{x_0} - 0.017{x_1} - 0.099{x_2},\\
1.43\lambda  \le 0.821 - 0.625{x_0} + 0.179{x_1} - 3.036{x_2},\\
1.143\lambda  \le 0.737 - 0.207{x_0} + 0.116{x_1} - 0.066{x_2},\\
\lambda  + d_{11}^ -  = 1,\,\,\,\lambda  + d_{12}^ -  = 1,\,\,\,\lambda  + d_{21}^ -  = 1,\,\,\,\\
\lambda  + d_{22}^ -  = 1,\,\,\lambda  + d_{31}^ -  = 1,\,\,\,\lambda  + d_{32}^ -  = 1,\\
{x_0} + {x_1} + {x_2} \le 5,\,\,\,{x_0} + {x_1} - {x_2} \le 2,\\
{x_0} + {x_1} + {x_2} \ge 1,\,\,\, - {x_0} + {x_1} + {x_2} \le 1,\\
{x_0} - {x_1} + {x_2} \le 4,\,\,\,\,\,{x_0} + 2{x_2} \le 4,{x_0} = 1.25,\\
d_{11}^ - ,d_{12}^ - ,d_{21}^ - ,d_{22}^ - ,d_{31}^ - ,d_{32}^ -  \ge 0\\
{x_0},{x_1},{x_2} \ge 0,\\
\lambda  \ge 0
\end{array} \right.
\end{array} \right.\]

The above model is a linear programming model which can be easily solved by Maple 18.02 optimization toolbox. The candidate optimal solutions are obtained as
\[x_{}^{S} =\left(x_{0}^{F} =1.25,x_{1}^{S} =0.75,x_{2}^{S} =0\right)\] 

(Step 8): Let the upper level decision maker be satisfied by the obtained candidate solution $x_{}^{S} =\left(x_{0}^{F} =1.25,x_{1}^{S} =0.75,x_{2}^{S} =0\right)$. 

(Step 9): Therefore, objective values are
\[f_{11} =-0.48,f_{12} =0.33,f_{21} =-0.45,f_{22} =-1.02,f_{31} =-0.35,f_{32} =0.61,\] 
and membership values are 
\[\mu _{11} \left(f_{11} \right)=0.83,\mu _{12} \left(f_{12} \right)=0.72,\mu _{21} \left(f_{21} \right)=0.47,\mu _{22} \left(f_{22} \right)=1,\mu _{31} \left(f_{31} \right)=0.43,\mu _{32} \left(f_{3} \right)=0.52.\] 

Baky [34] obtained the solution of the above problem with the use of the interactive fuzzy goal programming approaches. Comparative results are given in the following Table 4.

\begin{center}
Table 4 Comparison of solutions by different methods
\begin{tabular}{|p{0.9in}|p{1.6in}|p{0.9in}|} \hline 
 & The suggested approach\newline $\left(1.25,0.75,0\right)$ & Baky's approach [34]\newline $\left(1,0,0\right)$ \\ \hline 
$\mu _{11} \left(f_{11} \right)$  & $0.83$ & $0.46$ \\ \hline 
$\mu _{12} \left(f_{12} \right)$ & $0.72$ & $0.76$ \\ \hline 
$\mu _{21} \left(f_{21} \right)$ & $0.47$ & $0.31$ \\ \hline 
$\mu _{22} \left(f_{22} \right)$ & $1$  & $1$  \\ \hline 
$\mu _{31} \left(f_{31} \right)$ & $0.43$ & $0.54$ \\ \hline 
$\mu _{32} \left(f_{32} \right)$ & $0.52$ & $0.52$ \\ \hline 
\end{tabular}
\end{center}
From Table 4, all of the sums of the upper level's membership values produced by the proposed procedure is greater than that produced by Baky's solutions in [34]. Moreover,  all of the sums of the upper and lower level's objective values produced by the proposed procedure is smaller than (for minimization type objective) that produced by Baky's solutions  in [34]. 

Furthermore, the above example solved by Toksari and Bilim [36]. They used a neccessary constraints (i.e., $-\left(d_{1}^{ij} x_{1} +d_{2}^{ij} x_{2} +...+d_{m}^{ij} x_{m} \right)+D_{ij}^{-} \le \beta _{}^{ij} ,$) which is obtained from Baky's  linear approach in [34]. In addition, their results are $\mu _{11} \left(f_{11} \right)=1$and $\mu _{12} \left(f_{12} \right)=0.17.$ For these reasons, their solutions for the upper level decision maker is not most satisfied.

The advantage of the offered solution procedure over the current methods [34, 36]  is that there is no limitation the weights attached to the fuzzy control  operator $\lambda$ in the constraints. The claim that all fuzzy linear program has an equivalent weighted linear goal program where the weights are limited as the reciprocals of the acceptable violation constants are not every time actual. 

The other advantage of the proposed solution procedure is that the solution of any bi-level fuzzy fractional goal programming problems could be achieved efficiently without any computational complexities.  Therefore, these solutions show that the proposed bi-level fuzzy goal programming approach based on Taylor series in this paper is possible and feasible.

\section{Conclusion}

In this paper, a fuzzy goal programming algorithm based on Taylor series is introduced to obtain the optimal solution of decentralized bi-level multiobjective fractional programming (DBLMOFP) problems with a single decision maker at the upper level and multiple decision makers at the lower level. The proposed fuzzy goal programming solution procedure provides the most satisfactory solution for all the decision makers at both the levels by reaching the aspired levels of the membership goals. In this procedure, Taylor series approach is applied to linearize fractional membership functions associated with each objective. Also, the proposed solution procedure has an interactive structure as it provides the upper level to provide the opportunity for exchange the data presented that the upper level decision maker is not satisfied from this solution. Finally, application of the proposed solution procedure is handled with a known numerical example and then the effectiveness of the solutions shown by the proposed method is proved. Further, from table 4, the proposed solution procedure in this paper provides a more efficient solution comparing to the solution procedure of Baky [34].

\end{document}